\documentclass[10pt]{amsart}
\usepackage[cp866]{inputenc}
\usepackage[english]{babel}
\usepackage{amsmath}
\usepackage{amssymb}
\usepackage{amsfonts}

\newtheorem{lemma}{Lemma}
\newtheorem{theorem}{Theorem}

\newtheorem{proposition}{Proposition}

\newcommand {\bproof }{{\par\medskip\noindent \bf Proof. }}
\newcommand {\eproof }{\hfill $\blacktriangle$ \\ \medskip}

\title[  Propelinear perfect codes from
 the regular subgroups of $GA(r,q)$]{ q-ary propelinear perfect codes from
 the regular subgroups of
the $GA(r,q)$ and their ranks*}
\thanks{* This work was funded by the Russian Science Foundation under grant 22-21-00135,  https://rscf.ru/en/project/22-21-00135/}
\author{Ivan Mogilnykh}
\address{Sobolev Institute of Mathematics}
\email{ivmog84@gmail.com}
\date{}
\begin{document}

\maketitle

 \begin{abstract} We propose a new method of constructing $q$-ary propelinear perfect codes.
The approach utilizes permutations of the fixed length $q$-ary vectors that arise from the automorphisms
of the regular subgroups of the affine group. For any prime $q$ it is shown that the new class contains an infinite series of $q$-ary propelinear perfect codes of varying ranks of growing length. 
\end{abstract}


\section{Introduction}

Propelinear codes were introduced by Rif\`{a}, Basart and Huguet in \cite{Rifa1} and provide a general view on linear, additive (including $Z_4$-linear codes) and other
classes of codes. There are instances where propelinear approach yields codes with larger size than linear ones.
In particular, this holds for Preparata codes, and all known classes of these codes are shown to be propelinear \cite{HKCS}, \cite{BPRZ}, \cite{ZinZin}.   

Unlike binary codes, which are well-studied, there are rather few works devoted to $q$-ary propelinear codes, for $q\geq 3$. 
We refer to \cite{KP}, \cite{BMRS}, \cite{KP2}, \cite{ABE} for  $q$-ary perfect, MDS and generalized Hadamard propelinear codes.

In this work we propose a method of constructing $q$-ary, $q\geq 2$ perfect codes based on the automorphisms of the regular subgroups of
the general affine group $GA(r,q)$ and a Mollard construction \cite{Mol}.
 For $q=2$ this approach which uses particular case of Solov'eva construction was the topic of study of works \cite{MogSol} and \cite{MogSol2}. 
In \cite{MogSol} the values for the ranks and the kernels of these codes were found.  
A criteria for coordinate transitivity of resulting code in terms of double cosets of $GL(r,2)$ was suggested in \cite{MogSol2}. 
An infinite series of binary extended perfect codes  were constructed in \cite{MogSol2} with automorphism groups acting transitively on the code and transitively on the set of its neighbors.
Codes with such exquisite algebraic properties are known as neighbor-transitive \cite{GP}.

For a thorough study of ternary perfect codes of length $13$ we refer to \cite{KS}, where a large class of codes from construction \cite{Rom} were classified, and the values of their ranks and kernels were described. The rank problem for $q$-ary perfect codes was solved in \cite{PV}. 

Basic definitions are given in Section 2. In Section 3.1 the general  Mollard approach \cite{Mol} is described in terms of  work \cite{Rom}. The construction involves permutations of the vectors of $F_q^r$ .
In Section 3.2 we focus on the case when the permutations arise from the automorphisms of regular subgroups of $GA(r,q)$ and then the resulting perfect codes are propelinear. 
 Section 4 is on ranks of the codes from general Mollard construction, which could be described in terms of such permutations of $F_q^r$.
Section 5 and 6 are devoted to obtaining propelinear perfect codes of various ranks. The idea behind the last two Sections is a natural iterative approach for regular subgroups of $GA(r,q)$.

\section{Definitions}

 The all-zeros and all-ones vectors of the vector space $F_q^n$ are denoted as ${\bf 0}$, ${\bf 1}$ and their length will be clear
from the context. The concatenation of vectors $x$ and $y$ is denoted as $x|y$. For  $q$-ary codes $C$ and $D$, then we use the following notation $C
\times D$:
$$ C
\times D=\{(x|y):x \in C, y\in D \}.$$

 {\it The automorphism group} $Aut(F^n_q)$ of $F_q^n$ is defined as the group of all isomet\-ries of $F_q^n$, i.e. 
the automorphism group of the corresponding Hamming graph $H(n,q)$. {\it The automorphism group} $Aut(C)$ of a code $C$ is the setwise stabilizer of $C$ in $Aut(F^n_q)$.
 A $q$-ary {\it propelinear} code (original definition was given in \cite{Rifa1}) is a $q$-ary code whose automorphism group
contains a subgroup acting regularly on the codewords of $C$. An important invariant is the {\it rank} of a $q$-ary code $C$ which
is the dimension of its linear span over $F_q$. We denote the latter by $<C>$.

The {\it general affine group} $GA(r,q)$ is the group of all transformations $(a,M)$, where
$a\in F_q^r$, $M\in GL(r,q)$, acting on the column-vectors $b\in F_q^r$
as follows: $$(a,M)(b)=a+Mb,$$
with respect to the composition:
\begin{equation}
\label{GAcomp}
(a,M)(b,M')=(a+Mb,MM').\end{equation}

A subgroup of $GA(r,q)$ is called  {\it regular} if it acts regularly on the vectors of  $F_q^r$with respect to the
above defined action. Apart from the trivial example such as the translation group $(F_q,+)$ there 
are many regular subgroups of $GA(r,q)$. In Example 1 of this paper we consider a regular subgroup of $GA(2,q)$
for a  prime $q$, which is isomorphic to $(F_q,+)$ but no conjugate in $GA(2,q)$.

\section{The construction of q-ary propelinear perfect codes
from regular subgroups of $GA(r,q)$}
\subsection{Concatenation construction for $q$-ary perfect codes}
The construction for propelinear perfect codes is based on more a general method of Mollard \cite{Mol}.
We use the representation of this approach from work \cite{Rom} by Romanov.
\noindent Let $H_C$ be a parity check matrix of q-ary Hamming code C of length $\frac{q^r-1}{q-1}$, $H'$ be the matrix whose rows are all q-ary vectors of length $r$. A parity check matrix of $q$-ary Hamming code of length $\frac{q^{r+1}-1}{q-1}$ could be 
taken in a block form: \begin{equation}\label{mainr}\left(%
\begin{array}{c|c}
  {\bf 0} & {\bf 1} \\\hline
  H_C &  H' \end{array}%
\right).\end{equation}

For any $a\in F_q^r$  we use the notation below for a coset $C_a$ of the code $C$:

\begin{equation}\label{Ca}C_a=\{x:x \in F_q^{\frac{q^r-1}{q-1}}, H_C x^T=a\}.\end{equation}

We also denote by $D$ the linear code with the following parity check matrix: 
\begin{equation}\label{Hmat}H_D=\left(%
\begin{array}{c}
  {\bf 1} \\
 H' \\
\end{array}%
\right).\end{equation}
 We index the positions of $D$
with the columns of the parity check matrix $H_D$ and
the position has index $a$, $ a\in F_q^r$, if $\left(%
\begin{array}{c}
  1 \\
 a \\
\end{array}%
\right)$ is the corresponding column of $H_D$.
Denote by $e_a$ the vector of length $\frac{q^r-1}{q-1}$ of weight $1$ with one in the $a$th position.
For 
$a\in F_q^r$ denote by $D_a$ the coset $D+e_{\bf
0}-e_a.$
Note that for any $a$ the coset $D_a$ fulfills overall parity check, i.e. for any $y\in D_a$ we have ${\bf l}\cdot y^T={\bf 0}$.

\begin{theorem}\cite{Mol}\cite{Rom}\label{T1}
For any $q\geq 2$ and any permutation $\tau$ of the vectors of $F_q^r$ the code
$$S_{\tau}=\bigcup_{a\in F_q^r}C_{a}\times D_{\tau(a)}$$ is a $q$-ary
perfect code of length $\frac{q^{r+1}-1}{q-1}$.
\end{theorem}

Throughout the paper we assume that $\tau$ fixes the all-zero vector. 
 We note that the Hamming code with the parity check matrix (\ref{mainr}) is actually the code 
$\bigcup\limits_{a\in F_q^r}C_{a}\times D_{a}$.

\subsection{Concatenation construction for propelinear perfect codes from the regular subgroups of $GA(r,q)$}
Let $G$ be a regular subgroup of the general affine group $GA(r,q)$. Since $G$ acts
regularly on $F_q^r$, for any $a$ in $F_q^r$
 there is an element of $G$ that sends ${\bf 0}$ to $a$. We
denote this element by $g_a$. Since $g_a({\bf 0}) = a$ we see that the translation part of $g_a$ is $a$:
\begin{equation}\label{(4)}g_a = (a,M_a) \end{equation} 
for some nonsingular matrix $M_a$. Thus the elements of any regular subgroup of
$GA(r,q)$ are indexed by the vectors of $F_q^r$.

Let $T$ be an automorphism of $G$. The permutation $\tau$ of the vectors of $F_q^r$
such that  for any $a\in F_q^r$ $g_{\tau(a)} = T(g_a)$,
 is called the permutation {\it induced
by the automorphism} $T$. As any automorphism fixes the neutral element, the
permutation $\tau$ induced by any automorphism fulfills the equality $\tau({\bf 0}) = {\bf 0}$.

\begin{theorem} Let $\tau$ be the permutation induced by an automorphism of a regular subgroup of $GA(r,q)$. Then $S_{\tau}$ is a  $q$ ary propelinear perfect code of length $\frac{q^{r+1}-1}{q-1}$.
\end{theorem}

\section{The ranks of the codes obtained by concatenation construction}

Let $\tau$ be a permutation of the vectors of $F_q^r$ that fixes ${\bf 0}$.
 Since the positions of the code $D$ are indexed by the vectors of $F_q^r$,
 $\tau$ is also a permutation of the
positions of vectors of $D$. We define the distension of a permutation $\tau$ 
to be $dim(D)-dim(D \cap \tau(D))$. We have the following equality
\begin{equation}\label{(5)} dim(D)-dim(D \cap \tau(D))=rank\left(%
\begin{array}{c}
   H_D \\
   \tau(H_D) \\
\end{array}%
\right)-dim(D^{\perp}),\end{equation}
where $D^{\perp}$ is the dual code of $D$, i.e. the code whose generator matrix is the
parity check matrix $H_D$ for the code $D$.

Consider the code  $S_{\tau}=\bigcup\limits_{a\in F_q^r}C_{a}\times D_{\tau(a)},$ described in Section 3.1.  The main result of this section is the expression for the
rank of $S_{\tau}$ in terms of distention of $\tau$, which is given by exhibiting an explicit basis
of the linear span $<S_{\tau}>$.

For any $a\in  F_q^r$, we choose a representative of $C_a$ which we denote by $x_a$
throughout this section. The leader $e_{\bf 0} - e_a$ of the coset $D_a$ is denoted by $y_a$ and we have
$$ H_Dy_a^T=-\left(%
\begin{array}{c}
   0 \\
   a \\
\end{array}%
\right).$$

 In view of the
considered numeration of cosets of $C$ and $D$ via the vectors of $F_q^r$, we have the
following natural correspondence for linear dependencies in the coset spaces of
$C$ and $D$.

\begin{proposition}\label{stClDl} 
Let  $\tau$ be a permutation of $F_q^r$, fixing 
${\bf 0}$. For any elements $\alpha_a\in F_q$, $a\in F_q^r$ we have the following

$$ \sum_{a\in F_q^r}\alpha_a x_a\in C \mbox{ \, if and only if \,} \sum_{a\in F_q^r}\alpha_a y_{\tau(a)}\in\tau(D).$$

\end{proposition}
\bproof
The permutation $\tau$ fixes ${\bf 0}$ and acts on the positions of $F_q^r$
 that are
indexed by columns of $H_D$, i.e. the vectors of $F_q^r$. Therefore  we have:
$$\sum\limits_{a\in F_q^r}\alpha_a y_{\tau(a)}=\sum\limits_{a\in F_q^r}\alpha_a (e_{\bf 0}-e_{\tau(a)})=$$
$$\sum\limits_{a\in F_q^r}\alpha_a (e_{\bf \tau(0)}-e_{\tau(a)})=\tau(\sum_{a\in F_q^r}\alpha_ay_a).$$
 It remains to show that
$$\sum\limits_{a\in F_q^r}\alpha_a y_a\in D\mbox{ if and only if } \sum_{a\in F_q^r}\alpha_a x_a\in C.$$
Because the syndrome of $x_a\in C_a$ is $H_C x_a^T=a$, see (\ref{Ca}), we have that 
\begin{equation}\label{(6)}\sum_{a\in F_q^r} \alpha_a x_{a}\in C \mbox{ if and only if }H_C(\sum_{a\in F_q^r} \alpha_a x_{a})^T= \sum_{a\in F_q^r}\alpha_a a={\bf 0}.\end{equation}

Since $H_D(y_a)^T=-\left(%
\begin{array}{c}
   0 \\
 a\\
\end{array}%
\right)$, we obtain that
$$\sum_{a\in F_q^r} \alpha_a y_a\in D  \mbox{ holds if and only if }$$ $$ H_D(\sum_{a\in F_q^r} \alpha_a y_a)^T=-\left(%
\begin{array}{c}
   0 \\
 \sum_{a\in F_q^r}\alpha_a a \\
\end{array}%
\right)={\bf 0}.$$
This, combined with (\ref{(6)}), gives the required:
$$\sum_{a\in F_q^r} \alpha_a x_{a}\in C \Leftrightarrow  \sum_{a\in F_q^r}\alpha_a a={\bf 0}
\Leftrightarrow \sum_{a\in F_q^r} \alpha_a y_a\in D.$$

\eproof

In what follows we denote by  $z_1\ldots,z_{dim(C)}$  a basis of $C$, where $dim(C) =\frac{q^{r}-1}{q-1}-r$
 and by $v_1\ldots,v_{l}$ denote the vectors that complete a basis of $ D\cap \tau(D)$
to a basis of $D$. Note that here $l=dim(D)-dim(D)\cap dim(\tau(D))$ is the distension
of the permutation $\tau$.

We introduce three sets of vectors 
$$B=\{(x_a|y_{\tau(a)} ):a \in F_q^r\setminus {\bf 0}\},$$

$$B'=\{(z_i|{\bf 0}  ):i\in \{1,\ldots, dim(C)\}\},$$

$$B''=\{({\bf 0}|v_j   ):j\in \{1,\ldots,l\}\}.$$
We see that 
\begin{equation}\label{(7)}|B|+|B'|+|B''|=(q^r-1)+(\frac{q^r-1}{q-1}-r)+l=\frac{q^{r+1}-1}{q-1}-r-1+l.\end{equation}
We will now show that $B\cup B'\cup B''$ is a basis of the linear span  $S_{\tau}$.

\begin{lemma}\label{lemma 1}
The set $B\cup B'\cup B''$ is linearly independent.
\end{lemma}
\bproof
Clearly the sets $B\cup B'$ and $B''$ are linearly independent. Suppose that
$B\cup B'\cup B''$ is linearly dependent and consider a nonzero vector of the space
$<B\cup B'>\cap <B''>$. In view of $B$, $B'$ and $B''$ introduced above, the vector can be 
represented in two ways

$$\sum_{a \in F_q^r\setminus {\bf 0}} \alpha_a (x_a|y_{\tau(a)}   )+\sum_{i\in  \{1,\ldots, dim(C)\}} \beta_i(z_i|{\bf 0}  )=\sum_{j\in \{1,\ldots,l\}} \gamma_j({\bf 0}|v_j   )$$
for some  $\alpha_a,\beta_i,\gamma_j \in F_q$ and for all $i\in \{1,\ldots, dim(C)\}$, $j\in \{1,\ldots,l\}$, $a\in F_q^r$.
Equivalently, we have two equalities:

$$\sum_{a \in F_q^r\setminus {\bf 0}} \alpha_a x_a+\sum_{i\in  \{1,\ldots, dim(C)\}} \beta_i z_i={\bf 0},$$

$$\sum_{a \in F_q^r\setminus {\bf 0}} \alpha_ay_{\tau(a)}=\sum_{j\in \{1,\ldots,l\}} \gamma_j v_j   .$$

Since $z_i$'s are in $C$, the first of these equalities implies that
\begin{equation}\label{(8)}\sum_{a \in F_q^r\setminus {\bf 0}} \alpha_a x_a\in C.\end{equation}

By the choice of $\{v_j\}_{j\in \{1,\ldots,l\}}$ they complete a basis of $D\cap \tau(D)$ to a basis of $D$. Therefore their nontrivial linear combination $\sum\limits_{j\in \{1,\ldots,l\}} \gamma_j v_j$ is never in $\tau(D)$. Then the second equality gives that 

\begin{equation}\label{(9)}\sum_{a \in F_q^r\setminus {\bf 0}} \alpha_ay_{\tau(a)}\notin \tau(D).\end{equation}

Since (\ref{(8)}) and (\ref{(9)}) do not hold simultaneously by Proposition \ref{stClDl}, we obtain a contradiction.

\eproof
\begin{lemma}\label{lemma 2} Any vector of $B\cup B'\cup B''$ is in $S_{\tau}$ and $S_{\tau}\subseteq <B\cup B'\cup B''>$.
\end{lemma}
\bproof
Any vector of $B$ is $(x_a|y_{\tau(a)})$ for some $a\in F_q^r$ and therefore it is in
$S_{\tau}=\bigcup\limits_{a\in F_q^r} C_a\times D_{\tau(a)}$, so $B\subset S_{\tau}$. The set $B'$ is a basis $C \times 0$,
therefore it is included in $S_{\tau}$, whereas $B''$ completes a basis of ${\bf 0}\times(\tau(D) \cap D)$
to a basis of ${\bf 0}\times D$ and therefore $B''\subset C\times D\subset S_{\tau}$. Since $S_{\tau}$ is the union of
the cosets of $C\times D$ with representatives $(x_a|y_{\tau(a)} )$, $a\in F_q^r$, the said above
implies that it remains to prove that ${\bf 0}\times (\tau (D) \cap  D)$ is spanned by $B \cup B'$.

Given a vector $({\bf 0}|w)$, $w\in \tau (D) \cap D$ we will show that it is the sum of two
vectors from $<B>$ and $<B'>$. Recall that the parity check matrix $H_D$ of $D$ has an all-ones row, see (\ref{Hmat}). 
So, the code $D$, as well as $\tau(D)$, are subcodes of the supercode
with the parity check matrix $(1,\ldots,1)$. It is not hard to see that the vectors
$y_{\tau(a)}=e_{\bf 0} - e_{\tau(a)}$, for all $a\in F_q^r\setminus {\bf 0}$
 form a basis of the supercode. We consider a basis decomposition of $w\in \tau (D) \cap D$ on $y_{\tau(a)}$'s and
coefficients $\alpha_a\in F_q$, $a\in F_q^r$ such that
\begin{equation}\label{(10)} w=\sum_{a \in F_q^r\setminus {\bf 0}} \alpha_ay_{\tau(a)}.\end{equation}

Take the vector $\sum\limits_{a \in F_q^r\setminus {\bf 0}} \alpha_a (x_a|y_{\tau(a)})$, $w\in \tau(D)\cap D$ in $<B>$, which is
the linear combination of the vectors $(x_a|y_{\tau(a)})$, $a\in F^r_q\setminus {\bf 0}$ from $B$ with the
coeffcients $\alpha_a$'s. From the equality (\ref{(10)}) we see that the right side of this
vector is $w$:

\begin{equation}\label{(11)}\sum_{a \in F_q^r\setminus {\bf 0}}(\alpha_a x_a|\alpha_ay_{\tau(a)})=
(\sum_{a \in F_q^r\setminus {\bf 0}}\alpha_a x_a|w)\in <B>.\end{equation}

By the choice of the vector $w$, it is in $\tau(D)$. Moreover, $w$ decomposes as $\sum\limits_{a \in F_q^r\setminus {\bf 0}} \alpha_ay_{\tau(a)}$ and by Proposition
\ref{stClDl} the vector $\sum\limits_{a \in F_q^r\setminus {\bf 0}}\alpha_a x_a$ is in $C$. Because $B'$ is a basis of
$C\times 0$  we have that

$$\sum_{a \in F_q^r\setminus {\bf 0}}(\alpha_a x_a|{\bf 0})\in <B'>.$$

This, combined with (\ref{(11)}), gives that $({\bf 0}|w)$ is in $< B \cup B'>$. We conclude that
$C \times D$ and $S_\tau$ are in the span of $B\cup B'\cup B''$.

\eproof
From Lemmas \ref{lemma 1} and \ref{lemma 2} and equality (\ref{(7)}) we obtain the following.

\begin{theorem}\label{Trank}

Let  $\tau$ be a permutation of the vectors of $F_q^r$ with distension $l$ such that 
$\tau({\bf 0})={\bf 0}$. Then the rank of $S_{\tau}$ of length $\frac{q^{r+1}-1}{q-1}$ is equal to  $\frac{q^{r+1}-1}{q-1}-r-1+l$.

\end{theorem}

\section{The distension of the iteration of permutations}

Let $\tau_1$ and $\tau_2$ be permutations of $F^{r_1}_q$ and $F^{r_2}_q$ respectively, $\tau_1({\bf 0})={\bf
0}$, $\tau_2({\bf 0})={\bf 0}$. We represent any column-vector of $F_q^{r_1+r_2}$
as a concatenation $(^a_b)$ of some column-vectors $a\in F^{r_1}_q$ and $b\in F^{r_2}_q$. {\it The iteration of permutations} $\tau_1$ and $\tau_2$, denoted $\tau_1| \tau_2$ acts on the vectors of $F_q^{r_1+r_2}$ as follows:

 \begin{equation}\label{tausigma}(\tau_1| \tau_2)(^a_b)=(^{\tau_1(a)}_{\tau_2(b)}), \mbox{ for all } a\in F_q^{r_1}, b\in F_q^{r_2}.
\end{equation}

We show that the iterations of permutations that are induced by automorphi\-sms
of regular subgroups of $GA(r_1,q)$ and $GA(r_2,q)$ is a permutation induced by an
automorphism of a certain regular subgroup of $GA(r_1+r_2,q)$.
Let $G_1$ and $G_2$ be regular subgroups of $GA(r_1,q)$ and $GA(r_2,q)$. For elements $(a,M)\in G_1$ and $(b,M')\in G_2$ consider the following affine transformation from $GA(r_1+r_2,q)$ which we denote by $(a,M_1)\otimes (b,M_2)$:

$$((^a_b), \left(%
\begin{array}{cc}
  M_1& 0\\
      0& M_2\\
\end{array}%
\right)).$$

It is not hard to see that the direct product $\{(a,M_1)\otimes (b,M_2):(a,M_1)\in G_1,(b,M_2)\in G_2 \}$ of the groups $G_1$ and $G_2$ is a regular subgroup of $GA(r_1+r_2, q)$,
see e.g. \cite{Huffman}[Section 6]. Denote this group by $G_1 \otimes
 G_2$.

Let $T_1$ and $T_2$ be automorphisms of the groups $G_1$ and $G_2$ with induced
permutations $\tau_1$ and $\tau_2$ respectively. We define the permutation $T_1\otimes  
 T_2$ on the
elements of $G_1\otimes
G_2$ as follows: $(T_1\otimes
T_2)(g_1\otimes
g_2) = T_1(g_1)\otimes
T_2(g_2)$. It is obvious
that $T_1 
\otimes T_2$ is an automorphism of the group $G_1 
 \otimes G_2$ and the permutation
$T_1 \otimes 
 T_2$ of $F_q^{r_1+r_2}$ is $\tau_1| \tau_2$, defined earlier in (\ref{tausigma}). Thus we obtain the following.

\begin{proposition}\label{statt} Let $\tau_1$ and $\tau_2$ be permutations of $F_q^{r_1}$
 and $F_q^{r_2}$ induced by auto\-morphisms
of regular subgroups of $GA(r_1, q)$ and $GA(r_2, q)$, $q\geq 2$. Then $\tau_1|\tau_2$ is the permutation
induced by an automorphism of the regular subgroup $G_1$ 
 $G_2$ of $GA(r_1 + r_2, q)$
and the code $S_{\tau_1|\tau_2}$ is propelinear.\end{proposition}
We leave the proof of the following theorem without proof.

\begin{theorem}\label{Tsum} Let $\tau_1$ and $\tau_2$ be permutations of $F_q^{r_1}$
 and $F_q^{r_2}$, $q\geq 2$ with distensions
$l_1$ and $l_2$, respectively, such that $\tau_1({\bf 0}) = {\bf 0}$, $\tau_2({\bf 0}) = {\bf 0}$. Then the distension of the
permutation $\tau_1|\tau_2$ is $l_1 + l_2$.
\end{theorem}
\section{An infinite series of propelinear perfect codes with different
ranks}

We start this section with an example.

{\bf Example 1}. Let $q$ be a prime, $q\geq 3$. We will now construct a regular
subgroup of $GA(2,q)$ isomorphic to $Z_q^2$ but not conjugate to the translation
group $(F_q^2, +)$ in $GA(2, q)$. We then show that there is an automorphism of this
group with induced permutation of the vectors $F_q^r$
 having distension $2$.

Consider the following affine transformations
$$g=((^1_0),Id), h=((^0_1),\left(%
\begin{array}{cccccc}
  1 & 2 \\
  0 & 1 \\
\end{array}%
\right)).$$
It is not hard to see that $g$ and $h$ commute. Moreover, the
following holds:

\begin{equation}\label{eqgh}g^ih^j=((^i_0),Id)((^{j(j-1)}_j), \left(%
\begin{array}{cccccc}
  1 & 2i \\
  0 & 1 \\
\end{array}%
\right))=((^{i+j(j-1)}_j), \left(%
\begin{array}{cccccc}
  1 & 2i \\
  0 & 1 \\
\end{array}%
\right)).\end{equation}

For distinct pairs $(i,j)$ and $(i',j')$ the vectors $(^{i+j(j-1)}_j)$ and $(^{i'+j'(j'-1)}_{j'})$
are different. Therefore the group spanned by $g$ and $h$ is a regular subgroup of
$GA(2, q)$, isomorphic to $Z_q^2$.

Consider the permutation $T$ of the elements of the subgroup spanned by $g$
and $h$ defined as follows:
 $$T(g^ih^j)=g^jh^i$$
for all $i,j\in \{0,\ldots,q-1\}$. Since the $g$ and $h$ commute, the involution $T$ is an
automorphism of the group spanned by $g$ and $h$. By definition, the induced permutation $\tau$ of the automorphism $T$ is such that $\tau(a) = b$ if $T((a,M)) = (b,M')$ where $(a,M)$ and $(b,M')$ are elements of the considered
regular subgroup. From (\ref{eqgh}) we have $$g^ih^j=((^{i+j(j-1)}_j), \left(%
\begin{array}{cccccc}
  1 & 2i \\
  0 & 1 \\
\end{array}%
\right))$$ and because $g$ and $h$  commute, we obtain $$h^ig^j=((^{j+i(i-1)}_i), \left(%
\begin{array}{cccccc}
  1 & 2j \\
  0 & 1 \\
\end{array}%
\right))$$ for all $i,j\in \{0,\ldots,q-1\}$ and therefore we have  $$\tau(^{i+j(j-1)}_j))=(^{j+i(i-1)}_i).$$
In particular if pairs $(i,j)$ are equal to  $(1,0)$, $(0,1)$,
$(-1,-2)$ and $(0,2)$, we obtain \begin{equation}\label{eqtauex} \tau(^1_0)=(^0_1), \tau(^0_1)=(^1_0),
 \tau(^5_{-2})=(^0_{-1}),
\tau(^2_2)=(^2_0).\end{equation}

Using (\ref{(5)}) we find the distension of $\tau$, i.e.  $rank\left(%
\begin{array}{c}
   H_D \\
 \tau (H_D)\\
\end{array}%
\right)-dim(D^{\perp})=rank\left(%
\begin{array}{c}
   H_D \\
 \tau (H_D)\\
\end{array}%
\right)-3$. Since all-ones vectors are rows of both $H_D$ and
$\tau(H_D)$, we have that $$rank\left(%
\begin{array}{c}
   H_D \\
 \tau (H_D)\\
\end{array}%
\right)=1+rank \left(%
\begin{array}{cccc}{\bf 0}& b^2 & \ldots &  b^{q^{2}}\\{\bf 0}& \tau(b^2) & \ldots &  \tau(b^{q^{2}})\\
\end{array}%
\right),$$ where $b^2, \ldots,b^{q^2}$ are nonzero vectors of $F_q^2$. We take the first four nonzero $b_i$'s as follows:
$$b^2=(^1_0), b^3=(^0_1), b^4=(^5_{-2}), b^5=(^2_2).$$

From (\ref{eqtauex}) applying elementary transformations to the rows of the matrix we
see that $$rank \left(%
\begin{array}{cccc}{\bf 0}& b^2 & \ldots &  b^{q^{2}}\\{\bf 0}& \tau(b^2) & \ldots &  \tau(b^{q^{2}})\end{array}%
\right)=
rank \left(%
\begin{array}{ccccc} 1 & 0 & 5&2 &   \ldots \\
 0 & 1 & -2&2 &   \ldots \\
0 & 1 & 0&2 &   \ldots \\
1 & 0 & -1&0 &   \ldots \\
\end{array}%
\right)=$$ $$rank \left(%
\begin{array}{ccccc} 1 & 0 & 5&2 &   \ldots \\
 0 & 1 & -2&2 &   \ldots \\
0 & 0 & 2&0 &   \ldots \\
0 & 0 & 6&2 &   \ldots \\
\end{array}%
\right) =4 $$
and conclude that $\tau$ is of distension 2.
\begin{theorem} For all prime $q, q \geq 3$, $r\geq 2$ and $i\in \{0,\ldots, \lfloor r/2 \rfloor \}$ there is a
propelinear $q$-ary code $S_{\tau}$ of length $\frac{q^{r+1}-1}{q-1}$  and rank $\frac{q^{r+1}-1}{q-1}-r-1+2i$.
\end{theorem}

\bproof
Let $\tau$ be an induced permutation of $F_q^2$ with distension $2$ from Example 1. The permutation $\tau|\ldots |
\tau| id|\ldots | id$ of $F_q^r$, where $\tau$
is taken $i$ times, and identity is taken $r-2i$ times. From Proposition
\ref{statt} the code $S_{\tau|\ldots |
\tau| id|\ldots | id}$ is propelinear. In view of Theorem \ref{Tsum} 
the distension of $\tau|\ldots |
\tau| id|\ldots | id$ is $2i$, so from Theorem \ref{Trank} we obtain the desired value for rank.

\eproof

\end{document}